\documentclass[secthm]{elsart}
\usepackage{natbib}
\usepackage{amssymb}
\usepackage{upref} 
\usepackage[british]{babel}
\usepackage{amsmath,amscd,amsfonts}  
\usepackage[T1]{fontenc}
\usepackage[latin1]{inputenc}
\usepackage{latexsym}
\usepackage[dvips]{graphics}
\usepackage{enumerate}
\usepackage{url}
\usepackage{amsbsy}
\usepackage{xspace}

\journal{Statistics and Probability Letters}

\newcounter{myexamples}

\newenvironment{exemple}[1][]%
{\par {\bf %
\noindent #1 
} \nopagebreak  \\ \nopagebreak }{ }

\newenvironment{exempleNUM}[1][]%
{\par \addtocounter{myexamples}{1} {\bf %
\noindent Counter-Example~#1 {\Alph{myexamples}}
} \nopagebreak  \\ \nopagebreak }{ }

\newcommand{\ie}{{\it i.e. }}
\newcommand{\impt}[1]{{\em #1}}

\renewcommand{\leq}{\leqslant}
\renewcommand{\geq}{\geqslant}

\newcommand{\miniop}[3]{%
\renewcommand{\arraystretch}{0.6}
\begin{array}{c}
{\scriptstyle #1}\\
#2\\
{\scriptstyle #3}
\end{array}
\renewcommand{\arraystretch}{1}}

\newcommand{\card}[1]{\ensuremath{\# #1}}

\renewcommand{\L}{\Lambda}
\newcommand{\s}{\ensuremath{\sigma}}
\renewcommand{\e}{\ensuremath{\eta}}
\renewcommand{\r}{\ensuremath{\rho}}
\renewcommand{\o}{\ensuremath{\omega}}
\renewcommand{\t}{\ensuremath{\tau}}
\renewcommand{\d}{\ensuremath{\delta}}
\renewcommand{\b}{\ensuremath{\beta}}
\newcommand{\z}{\ensuremath{\zeta}}

\renewcommand{\epsilon}{\varepsilon}

\newcommand{\Zd}{{\Zset}^d}

\newcommand{\N}{\mathbb N}

\newcommand{\tq}{ : \ }

\newcommand{\ind}[2]{1 \hspace{-0.9ex} 1_{#1} {#2}}

\newcommand{\spin}{\ensuremath{S}\xspace}

\newcommand{\cvinf}{\ensuremath{\spin^{\Zd}}}

\newcommand{\cvinfLx}[1]{\ensuremath{\spin^{#1}}}

\newcommand{\SpinConfigSpace}{\ensuremath{\{-1,+1\}^{\Zd}}}

\newcommand{\abs}[1]{\ensuremath{ \Big|  #1  \Big| }}

\newcommand{\norme}[1]{\ensuremath{ \| #1 \|_{_{_1}} }}

\newcommand{\tvert}{\vert \hspace{-0.8ex} \parallel}

\newcommand{\ntrois}[1]{\tvert #1 \ \tvert}

\newcommand{\BL}{\ensuremath{ {{\mathcal B}({L})}}}

\newcommand{\nusup}{\ensuremath{{\nu^\maximalConf}}}
\newcommand{\nuinf}{\ensuremath{{\nu^\minimalConf}}}

\renewcommand{\preceq}{\preccurlyeq}
\renewcommand{\succeq}{\succcurlyeq}

\newcommand{\proj}[2]{{\wp}_{#2} \ {#1}}

\newcommand{\Q}{\ensuremath{\mathbf Q}}

\newcommand{\coup}[1]{\mathbf{P}^{#1}}
\newcommand{\coupling}[2]{\coup{} \Big( #1 \Big| #2 \Big)}

\newcommand{\opreceq}{\circledast}

\newcommand{\K}{\ensuremath{{\mathcal K}}}

\renewcommand{\S}{\ensuremath{{\mathcal S}} }

\newcommand{\be}[1]{
\begin{equation}\label{#1}}
\newcommand{\ee}{\end{equation}}

\newcommand{\minimal}{\ensuremath{{-}}}
\newcommand{\maximal}{\ensuremath{{+}}}
\newcommand{\minimalConf}{\ensuremath{{\pmb{-}}}}
\newcommand{\maximalConf}{\ensuremath{{\pmb{+}}}}

\newcommand{\Var}[2]{\ensuremath{\Delta_{#2}(#1)}}

\newcommand{\semi}[1]{\ensuremath{(\leftarrow,#1]}}

\newcommand{\Aposet}{\ensuremath{\mathcal A}}

\newcommand{\Sposet}{\ensuremath{\mathcal S}}

\newcommand{\Spreceq}{\ensuremath{{\preceq}_{\mathcal S} \ }}

\newcommand{\Apreceq}{\ensuremath{{\preceq}_{\mathcal A} \ }}

\newcommand{\Pcourt}{\ensuremath{\miniop{}{\otimes}{k \in \Zd} p_k}}

\newcommand{\preceqlino}{\ensuremath{\leq_{n}}}

\begin{document}

\footnotesize

\begin{frontmatter}

\title{Increasing coupling of \\ Probabilistic Cellular Automata}

\date{2005}

\author{Pierre-Yves {Louis}}
\ead{louis@math.uni-potsdam.de}
\ead[url]{http://www.math.uni-potsdam.de/$\sim$louis}

\address{Institut für Mathematik, Potsdam Universität, \\ Am neuen Palais, Sans Souci, %
\\ Potsfach 60 15 53, D-14 415 Potsdam}

 Post-print\footnote{\url{http://www.sciencedirect.com/science/article/pii/S0167715205001501?np=y}}
 of Statistics \& Probability Letters \\
 74 (1), 1 August 2005, pp. 1--13, ISSN 0167-7152 \\
doi:10.1016/j.spl.2005.04.021\\
{\footnotesize Creative Commons Attribution Non-Commercial No Derivatives License}

\begin{abstract}
\begin{minipage}{17cm}
We give a necessary and sufficient condition for the
existence of an increasing  coupling of $N$ ($N \geq 2$)
 synchronous dynamics on
\cvinf (PCA).
Increasing means the coupling preserves stochastic ordering.
We first present our main construction theorem in the case where
\spin is totally ordered; applications to attractive PCA's are given.
When \spin is only partially ordered,
we show on two examples that a coupling of more than two synchronous dynamics may not exist.
We also prove an extension of our main result for a particular class of partially ordered spaces.
\end{minipage}
\end{abstract}


\end{frontmatter}

%
%

\section{Introduction}
Probabilistic Cellular Automata (abbreviated in PCA) are discrete-time Markov
chains on a product space $S^{\L}$ (\impt{configuration space})
whose transition probability is a \impt{product measure}.
$\spin$ is assumed to be a finite %
set (\impt{spin space}). We denote by $\L$
(set of \impt{sites}) a subset, finite or infinite, of $\Zd$.
Since the transition probability kernel $P(d\s|\s')$ %
($\s,\s' \in
S^{\L}$) is a product measure, all interacting elementary components (\impt{spins}) %
 $\{\s_{k}: k
\in \L\}$ are simultaneously and independently updated %
(\impt{parallel updating}). This synchronous transition is the main feature of PCA and differs from the
one in the most common Gibbs samplers, where only
one site is updated at each time step (\impt{sequential updating}).
In opposition to these sequential updating dynamics, it is simple to
define PCA's on the infinite set $\cvinf$ without passing to
continuous time. \\ \indent
Probabilistic Cellular Automata were first studied as Markov chains
in the 70's (see~\cite{ToomEtAl}).
We refer for instance to~\cite{ThesePyl} for a %
recent historical overview and a list of
applications of Cellular Automata dynamics, which 
are to be found in physics, biology, image restoration (see \cite{YounesImage})...
PCA dynamics may present a variety of behaviours.
Let us only mention the following: contrarily to the usual
discrete time sequential updating dynamics, for a given measure $\mu$, %
there is no canonical way of constructing a PCA for which $\mu$ is stationary. Moreover, there exist %
Gibbs measures on $\cvinfLx{\mathbb Z^2}$ such that no PCA admits them as stationary reversible measures %
(see Theorem~4.2 in \cite{Dawson74}). \\ \indent
Coupling refers to the
construction of a product probability space on which several dynamics may evolve simultaneously, and having the %
property that the marginals coincide with each one of these dynamics.
Coupling techniques for stochastic processes are now well established, powerful tools %
of investigation. 
We refer to \cite{LindvallBook} and \cite{Thorisson} for a more extensive
review and applications to a large scope of probabilistic objects.
The first use of a coupling of Probabilistic Cellular Automata is to be found in~\cite{V69}. %
It was also used in~\cite{MaesRMP93}.
Recently, the coupling constructed in this paper was used to state some necessary and sufficient condition
for the exponential ergodicity of attractive PCA's (see \cite{PylExpErg}).
This last result relies on the fact that our coupling preserves a stochastic order between the configurations
(so called \impt{increasing coupling}).
In \cite{LopezSanz2000}, the authors gave necessary and sufficient condition for the existence of
a  coupling preserving   the stochastic order
between two PCA's on \cvinf , where $\spin$ is a partially ordered set.
In this paper, we give a necessary and sufficient condition for the existence of an increasing coupling
of any finite number of possibly different PCA dynamics.
As some counter examples will show,
there is a gap between the construction of an increasing coupling of two PCA's and that of an increasing coupling of~$N$
PCA's with $N \geq 3$.
Moreover, we give here an explicit algorithmic construction of this coupling, which
is a kind of graphical construction.
We also give several examples and general applications of the constructed coupling.
Indeed, the motivation for coupling together three or more PCA's comes from the paper~\cite{PylExpErg}, where
a comparison between four different PCA dynamics proved to be useful. \\ \indent
In section~\ref{section_main} we state our main result, namely
the existence, under some necessary and sufficient condition of monotonicity %
(Definition~\ref{DefIncreasingNUple}), of an
increasing coupling of several PCA dynamics (Theorem~\ref{ExistenceCouplageCroissant}).
Corollary~\ref{Corollary_main} states the existence of some universal coupling of any attractive PCA and some
significant examples are also presented.
In section~\ref{proofs_main} we prove these results, and state some important %
property (Lemma~\ref{Propriete_Compatibilite}) of coherence between the different couplings.
In section~\ref{Applications}, we then present some useful applications of the coupling just constructed.
In section~\ref{section_ordre_partiel} we consider the case where $\spin$ is a partially ordered set.
Two counter-examples show that it may happen that an increasing
coupling of $N$ PCA dynamics does not exist when $N \geq 3$. %
A generalisation of Theorem~\ref{ExistenceCouplageCroissant}
and Corollary~\ref{Corollary_main} to the case where \spin is partially linearly ordered is presented. \\ \indent
Finally, let us point out that our motivation for considering partially ordered spin spaces comes from the study of
'block dynamics', where the siteds are not updated individually, but rather blockwise.
This amounts to consider PCA's on $(\spin^r)^{\Zd}$ where $r$ is the number of sites
in these blocks (even if \spin is totally ordered, $\spin ^r$ is not totally ordered in a natural way).

\newpage
\section{Definitions and main results\label{section_main}}

Let $\spin$ be a finite set, with a partial order
denoted by $\preceq$. The conjunction of $s \neq s'$ and $s \preceq s'$ will be denoted by $s \precneqq s'$.
Let $P$ denote a PCA dynamics on the product space
$\cvinf$, which means a time-homogeneous Markov
Chain on $\cvinf$ whose transition probability kernel $P$ verifies,
 for all configurations $\eta \in
\cvinf$, \mbox{$\s=(\s_k)_{k \in \Zd} \in \cvinf$},
$P(\ d\s \ | \ \e \ )=\miniop{}{\otimes}{k \in \Zd} p_k(\ d\s_k \ | %
\ \e \ )$,
where for all site $k \in \Zd$, $p_k(\ .\ |\e)$ is a probability measure on $\spin$, called
\impt{updating rule}.
In other words, \impt{given the previous time step} $(n-1)$,
all the spin values $(\o_k(n))_{k \in \Zd}$ at time~$n$
 are
\impt{simultaneously and independently updated}, each one according to
  the probabilistic rule $p_k(\ .\ |\ (\o_{k'}(n-1))_{k' \in \Zd})$.
  We let  $P=\Pcourt$.
  All PCA dynamics considered in this paper are
local, which means
$\forall k \in \Zd, \exists \ V_k \Subset \Zd, %
p_k(\ .\ |\e)=p_k(\ .\ |\e_{V_k})$.
The notation \mbox{$\L \Subset\Zd$} means $\L$ is
a finite subset of $\Zd$.
For any subset $\Delta$ of $\Zd$ and for all configurations
$\s$ and $\e$ of $\cvinf$, the configuration
$\s_\Delta \e_{\Delta^c}$ is defined by $\s_k$ for $k \in \Delta$,
$\e_k$ elsewhere.
We also let $\s_\Delta := (\s_k)_{k \in \Delta}$ too. \\ \indent
All the measures considered in this paper are probability measures.
For a probability measure~$\nu$ on~$\cvinf$
(equipped with the Borel $\s$-field associated to the product
topology),
$\nu P$ refers to the law at time~$1$ of the PCA dynamics with
law~$\nu$ at time~$0$:
$\nu P (d\s)=\int P(d\s |  \e ) \nu(d\e)$.
Recursively, $\nu P^{(n)}=(\nu %
P^{(n-1)}) P$ is the law at time $n$ of the system evolving
according to the %
PCA dynamics $P$ and having initial law $\nu$.
For each measurable function~$f: \cvinf \rightarrow \mathbb R_+$, $P(f)$ denotes the function on $\cvinf$
defined by $P(f)(\e)=\int f(\s) P(d\s|\e)$. \\ \indent
Let us now define basic notions of stochastic ordering $\preceq$.
Two configurations $\s$ and $\e$ of $\cvinfLx{\L}$ (with $\L \subset \Zd$)  %
satisfy $\s \preceq \e$ if $\forall k \in \L, \s_k \preceq \e_k$.
A real function $f$ on $\cvinfLx{\L}$ will be increasing if $\s \preceq \e \Rightarrow $ %
$f(\s) \leq f(\e)$. Thus two probability measures $\nu_1$ and $\nu_2$ satisfy the %
stochastic ordering $\nu_1 \preceq \nu_2$ if,
for all increasing functions~$f$ on $\cvinfLx{\L}$, $\nu_1(f) \leq \nu_2(f)$, with the %
notation $\nu(f)=\int f(\s) \nu (d\s)$.
Considered as a Markov chain, a PCA dynamics $P$ on $\cvinfLx{\L}$ ($\L \subset \Zd$) is said to be %
\impt{attractive} if
for all increasing functions $f$, $P(f)$ is still increasing. This requirement is equivalent to
($\mu_1 \preceq \mu_2 \Rightarrow  \mu_1 P \preceq \mu_2 P$), where $\mu_1,\mu_2$ are two %
probability measures on $\cvinf$.
\begin{defn}[Synchronous coupling of PCA dynamics] \hspace{2mm} \\ \indent
Let $P^1,P^2,\hdots,P^N$ be $N$ probabilistic cellular automata dynamics,
with $P^i=\miniop{}{\otimes}{k \in \Zd} p^i_k$.
A \impt{synchronous coupling} of $(P^i)_{1\leq i \leq N}$  %
is a Markovian dynamics~$Q$ on~$(\cvinf)^N$,
which is also a PCA dynamics
whose  marginals coinncide respectively with $P^1,P^2,\hdots,P^N$.
Thus, $Q$~is such that $Q=\miniop{}{\otimes}{k \in \Zd} q_k$
and \\
$\forall i \in \{1,\hdots,N\}, \quad \forall s^i \in \spin, %
\ \forall \z^1,\hdots,\z^N \in \cvinf, $
\begin{equation}
 p^i_k(s^i \ | \ \z^i) = \sum_{s^j \in \spin, j \neq i} q_k\big( \ (s^1,\hdots,s^N) %
\ \big| \ (\z^1,\hdots,\z^N) \big).  \label{Marginales} %
\end{equation}
\end{defn}
\begin{defn}[Increasing $N$-tuple of PCA dynamics]
\label{DefIncreasingNUple}
Let $(P^1,P^2,\hdots,P^N)$ be an $N$-tuple of PCA dynamics,  %
where $N\geq 2$ and $P^i=\miniop{}{\otimes}{k \in \Zd} p^i_k$ ($1 \leq i \leq N $).
This $N$-tuple is said to be %
\impt{increasing} if
\begin{equation}
\z^1 \preceq \z^2 \preceq \hdots \preceq \z^N \quad \Rightarrow \quad %
 \label{Monotonie_globale}
P^1(\ . \ | \ \z^1 ) \preceq P^2(\ . \ | \ \z^2 ) \preceq \hdots %
\preceq P^N(\ .\  | \ \z^N ).
\end{equation}
\end{defn}
Since $P(\ . \ | \s)$ is a product measure, according to Proposition~2.9 in~\cite{ToomEtAl},
condition~(\ref{Monotonie_globale})
is equivalent
to: $\forall k \in \Zd,$
\begin{equation}
\z^1 \preceq \z^2 \preceq \hdots \preceq \z^N \quad \Rightarrow \quad %
\label{Monotonie_locale}
p^1_k(\ . \ | \ \z^1 ) \preceq p^2_k(\ . \ | \ \z^2 ) \preceq \hdots %
\preceq p^N_k(\ . \ | \ \z^N ).
\end{equation}
\begin{defn}[Increasing synchronous coupling]
A synchronous coupling~$\Q$ of an $N$-tuple $(P^i)_{1\leq i \leq N}$ of PCA dynamics
is said to be an
\impt{increasing coupling} of %
$(P^1, P^2, \hdots , P^N)$
if the following property holds:
for any initial configurations
 $ \s^1 \preceq \s^2 \preceq \hdots \preceq \s^N $, %
 for any time $n \geq 1 $,
\be{PreservationOrdreStochastique}
 \Q \ \big( \ %
\o^1(n) \preceq \hdots \preceq \o^N(n) \ \big| \ (\o^1,\hdots,\o^N)(0)= %
(\s^1,\hdots,\s^N) \ \big) =1  .
\ee
\end{defn}
We now state:
\begin{thm} 
\label{ExistenceCouplageCroissant} 
Let \spin \ be a totally ordered space. Let $(P^i)_{1\leq i \leq N}$
be an $N$-tuple of %
PCA dynamics on $\cvinf$.
There exists a synchronous coupling~$\Q$ of $(P^i)_{1\leq i \leq N}$
if and only if $(P^1,\hdots,P^N)$ is increasing.
\end{thm}
The increasing coupling we are constructing will be denoted by
\mbox{$P^1 \opreceq P^2 \opreceq \hdots \opreceq P^N$}.
Note that the property of preserving the order implies that the coupling has
the \impt{coalescence property}.
This means that if two components are taking the same value at some time, then they
will remain equal from this time on (as well as all the components in between).
In Lemma~\ref{Couplage_du_meme} we will see that
 if a PCA dynamics $P$ is attractive, then
for all $N \geq 2$, the $N$-tuple $(P,P,\hdots,P)$ is increasing.
As an immediate consequence of Theorem~\ref{ExistenceCouplageCroissant} and
Lemma~\ref{Couplage_du_meme} we then have the
\begin{cor} \label{Corollary_main} Let \spin \ be %
a totally ordered space, $P$ be a PCA dynamics on $\cvinf$ and $N \geq 2$.
There exists an increasing coupling
$P^{\opreceq N}$ of $(P,\hdots,P)$
if and only if $P$ is attractive.
\end{cor}
\indent Lemma~\ref{Fonction_repart_ordre_total} from section~\ref{proofs_main} gives %
a practical constructive criterion for testing
if an $N$-tuple of PCA dynamics is increasing or if a PCA is attractive.
We use it in the following examples.

\begin{exemple}[A family of PCA dynamics]
\indent Let $(P^{\beta_i,h_i})_{1 \leq i \leq N}$ be a family of
$N$ PCA dynamics on~$\SpinConfigSpace$, %
defined by
\mbox{$\forall k \in \Zd$}, \mbox{$\forall \e \in \{-1,+1\}^{\Zd}$},  %
\mbox{$\forall s \in \spin=\{-1,+1\}$}
\be{FormeLocaleNosPCAPrimaire}
 p^i_k(s \ | \ \e )=\frac{1}{2} \Big(1+ s \tanh ( %
\b_i \sum_{k'\in V_0} \K(k'+k) \e_{k'} + \b_i h_i  ) \Big),
\ee
where~$(\b_i)_{1 \leq i \leq N}$ are positive real numbers, $(h_i)_{1 \leq i \leq N}$ real numbers, $V_0 \Subset \Zd$ %
and
\mbox{$\K\ : \ V_0 \rightarrow \mathbb R$}
is an interaction function between sites
which is symmetric.

This example is important, since any reversible PCA dynamics on $\SpinConfigSpace$ can be presented in this form\footnote{ %
A PCA dynamics is said to be reversible if it admits
at least one reversible probability measure
(see subsection~4.1.1 in~\cite{ThesePyl}).}.
When $\b$ is fixed and $h_1 \leq \hdots \leq h_N$, the $N$-tuple $(P^{\beta,h_i})_{1 \leq i \leq N}$ is increasing.
On the other hand, note that in the case $h_i=0$, the assumption $\b_1 \leq \hdots \leq \b_N$ does not imply that
the $N$-tuple $(P^{\beta_i,0})_{1 \leq i \leq N}$ is increasing. Consider for instance $\b_1=\frac{1}{2}$,
$\b_2=3$, $d=2$, $\K$ such that $V_0=\{-e_1, e_1, -e_2, e_2\}$ where $(e_1,e_2)$ is a basis of $\mathbb R^2$.
Condition~(\ref{Monotonie_locale}) is false considering $k=0$, $\z^1_{V_0}$ consisting of four $-1$, and
$\z^2_{V_0}$ of three $-1$ and one $+1$.
\end{exemple}
\begin{exemple}[Example of an attractive PCA dynamics]
\indent Let $P^{\beta,h}$ be some PCA dynamics defined by the updating rule
(\ref{FormeLocaleNosPCAPrimaire}) ($\beta \geq 0$, $h \in \mathbb R$). We know from
Proposition~4.1.2 in \cite{ThesePyl} that
this dynamics is attractive if and only if $\K(.) \geq 0$.
For a more systematic study of this class, we refer to~\cite{DPLR} and \cite{PylExpErg}.
From Corollary~\ref{Corollary_main}, we can then construct an increasing coupling of such PCA.
We show in section~\ref{Applications} how this can be used.
\end{exemple}
\begin{exemple}[Example of an attractive PCA dynamics with $\mathbf{\card{\spin}=q}$, \mbox{$\mathbf{q \geq 2}$}]
\indent Let $\spin=\{ 1,\hdots,q \}$ ($q \geq 2$), and consider the updating rule \\
$\forall k \in \Zd, \quad \forall s \in \spin, \quad \forall \s \in \cvinf$,
$$ p_k(s|\s)= \frac{e^{\beta N_k(s,\s)}}{\sum_{s' \in \spin} e^{\beta N_k(s',\s)}}$$
where $\b \geq 0$, $V_k$ is a finite neighbourhood of~$k$ %
and $N_k(s,\s)$ is the number of $\s_{k'}$ ($k' \in V_k$) which are larger than~$s$.
This dynamics is attractive for any $\b$ non-negative.
\end{exemple}
%
%

\section{Proof of the main results\label{proofs_main}}
Assume in this section that $\spin$ is a totally ordered set.
Let us then enumerate the
spin set elements as
$ \spin = \{ \minimal, \hdots, s,s+1,\hdots, \maximal \}$,
where we denote with $\maximal$ (resp. $\minimal$)
the (necessarily unique) maximum (resp. minimum) value of $\spin$ and
for $s \in \spin$, $(s+1)$ denotes the unique element in $\spin$ such that
there is no $s'' \in \spin$, $s \precneqq s'' \precneqq s+1$.
A real valued function~$f$  on \cvinf is said to be \impt{local}
if
$ \exists \L_f \Subset \Zd,\ \forall \s \in \cvinf$,
\mbox{$f(\s)=f(\s_{\L_f}) $}.
\begin{lem} \label{Fonction_repart_ordre_total}
When \spin \ is a totally ordered space, the condition~(\ref{Monotonie_locale}) of monotonicity
is equivalent to \\
$  \forall k \in \Zd, \forall 
\z^1 \preceq \z^2 \preceq \hdots \preceq \z^N \in (\cvinf)^N, \forall s \in \spin $
\begin{equation} \label{ordre_fonction_repart}
F_k^1(s , \z^1 ) \geq F_k^2(s, \z^2 ) \geq \hdots %
\geq F_k^N(s ,\z^N ),
\end{equation}
where $F^i_k(s,\s)$ is
the distribution function  of $p^i_k(.|\s)$: %
$F^i_k(s,\s)= \miniop{}{\sum}{s'\leq s} p^i_k(s'|\s)$.
\end{lem}
\begin{pf}
The implication (\ref{Monotonie_locale}) $\Rightarrow$ (\ref{ordre_fonction_repart}) is straightforward
using the increasing function $f(s')=\ind{\{s'>s\}}{}$.
To prove (\ref{ordre_fonction_repart}) $\Rightarrow$ (\ref{Monotonie_locale}),
it is enough to remark that, for any function $f:\spin \to \mathbb R$,
\begin{equation} \label{Fonction_croissante_Fonction_repart}
 p^i_k(f|\s)= f(\maximal) + \miniop{}{\sum}{s \precneqq \maximal} (f(s) -f(s+1)) F^i_k(s,\s). \qed
\end{equation}
\end{pf}
\begin{pf*}{Proof of Theorem~\ref{ExistenceCouplageCroissant}} \\
Let us explain how to construct explicitly the increasing coupling
\mbox{$P^1 \opreceq P^2 \opreceq \hdots \opreceq P^N$}. \\
$n$ being a fixed time index, we need to describe the stochastic transition from
\mbox{$(\o^1,\hdots,\o^N)(n)$} (element of $\spin^N$) to
\mbox{$(\o^1,\hdots,\o^N)(n+1)$}.
Let $(U_k)_{k \in \L}$ be a family of independent random variables, %
distributed uniformaly on $]0,1[$.
Since we are  constructing  a synchronous coupling, it is enough %
to define the rule for a fixed site $k \in \Zd$.
Let~$r$ denote a fixed realisation of the random variable $U_k$ and
use the following \impt{algorithmic rule} to choose
the value $\o^i_k(n+1)$ for any $i \ (1\leq i \leq N)$:
\begin{equation} \label{RuleAlgorithm}
\left\{ \begin{array}{l}
\textrm{if } F^i_k (s-1,\o^i(n)) < r \leq  F^i_k(s,\o^i(n)) , \quad \minimal \precneqq s, \textrm{ assign }
\o^i_k(n+1)=s \\
\textrm{if } 0 \leq r \leq  F^i_k(\minimal,\o^i(n)) \textrm{ assign }
\o^i_k(n+1)=\minimal \ . \\
\end{array} \right.
\end{equation}
This rule corresponds to the definition of the coupling between times $n$ and $n+1$ according
to
\begin{equation} \label{regle_locale_algo}
\forall k \in \Zd, \quad  \Big( \o_k^i(n+1) \Big)_{1 \leq i \leq N} = %
 \Big( \big( F^i_k(\ . \ , \o^i(n)) \big)^{-1} (U_k) \Big)_{1 \leq i \leq N}
 \end{equation}
where $(F_k^i)^{-1}$ denotes the Lévy probability %
transform (generalised inverse probability transform) of the $F_k^i$ distribution function
$$ (F_k^i)^{-1}(t) = \miniop{}{\inf}{\preceq} \{ s \in \spin : F_k^i(s) \geq t \}, \quad t \in ]0,1[, \ i \in \{1,\hdots,N\}. $$
\indent Finally, remark that the stochastic dependence between
the components~$1\leq i \leq N$ comes from the fact that we are using
the \impt{same} realisation $r$ of $U_k$ for \impt{all} components.
It is easy to check that this coupling
preserves stochastic ordering assuming that $(P^1,\hdots,P^N)$ is increasing, since it is equivalent %
to check~(\ref{ordre_fonction_repart}) (Lemma~\ref{Fonction_repart_ordre_total}). \\
\indent Conversely, the condition~(\ref{ordre_fonction_repart}) is necessary. Assume the existence of a synchronous coupling
$(q_k)_{k \in \Zd}$ of $N$ PCA dynamics on $\cvinf$ which preserves stochastic ordering.
This means that for $\z^1 \preceq \hdots \preceq \z^N$, $q_k(\ . \ | (\z^1,\hdots,\z^N)) >0$ only on $(\spin^N)^+$,
where $(\spin^N)^+$ is the subset $\{(s^1,\hdots,s^N): s^1 \preceq \hdots \preceq s^N \}$ of $\spin^N$.
Let $s\in \spin$, $1\leq i <N$, and $\z^1 \preceq \hdots \preceq \z^N$ be fixed.
Using the condition~(\ref{Marginales}) on the $i$-th marginal of a coupling, we have
$$ F^i_k(s,\z^i)=\miniop{}{\sum}{(s^1,\hdots,s^N) \in A^i_s } q_k( (s^1,\hdots,s^N) | (\z^1,\hdots,\z^N)) ,$$
where $A^i_s=\{(s^1,\hdots,s^N) \in (\spin^N)^+: s^i \preceq s \}$.
Decompose
$A^i_s=A^{i+1}_s \sqcup \Delta^i_s$ with
\mbox{$\Delta^i_s=\{(s^1,\hdots,s^N) \in (\spin^N)^+:$} \mbox{$s^i \precneqq  s \preceq s^{i+1} \}$}
($\sqcup$ denotes the disjoint union).
Finally note that
$$ F_k^{i} (s,\z^{i}) =  F_k^{i+1} (s,\z^{i+1} ) + \miniop{}{\sum}{(s^1,\hdots,s^N) \in \Delta^{i}_s } %
q_k( (s^1,\hdots,s^N) | (\z^1,\hdots,\z^N)) $$
where the last term is non-negative. \qed
\end{pf*}
\indent One may readily notice the
\impt{compatibility property} satisfied by this coupling, and whose proof
is straightforward according to its construction.
\begin{lem} \label{Propriete_Compatibilite}
Let $N$ and $N'$ be two integers such that
$1 \leq N < N' $. Let
 $(P^1,\hdots,P^{N'})$ be $N'$  PCA dynamics.
The projection of the coupling $P^{1} \opreceq P^{2} %
\hdots \opreceq P^{N'}$ on any $N$ components $(i_1,\hdots,i_N)$ coincides with the
coupling $P^{i_1} \opreceq \hdots \opreceq P^{i_N}$.
\end{lem}
\begin{lem} \label{Couplage_du_meme}
Let $P$ be a PCA dynamics on \cvinf. It is an attractive dynamics
if and only if, for all $N \geq 2$, the $N$-tuple $(P,P,\hdots,P)$ is increasing.
\end{lem}
\begin{pf}
Assume $P$ is attractive. Let $k \in \Zd$ be fixed, and let $f_0$ be an increasing function on~\spin.
We consider the function $f$ on~\cvinf defined by $f(\s)=f_0(\s_k)$, $\forall \s \in \cvinf$.
Since $P(f)=p_k(f_0)$ is an increasing function, relation~(\ref{Monotonie_locale}) holds with $p_k^i=p_k, \forall i$.
The equivalence   (\ref{Monotonie_locale})~$\iff$~(\ref{Monotonie_globale}) %
gives $(P,\hdots,P)$ increasing for any~$N \geq 2$. \\ \indent
Conversely, assume $(P,P)$ is increasing. Then relation~(\ref{ordre_fonction_repart}) holds with %
the same dynamics on the two components. Let $f$ be an increasing function on~\cvinf such that $\exists k \in\Zd, %
\ \forall \s \in \cvinf$,
$f(\s)=f(\s_k)$. According to formula~(\ref{Fonction_croissante_Fonction_repart}), we conclude
that $P(f)$ is increasing. Recursively, we can state the same result for all local functions,
because of the product form of the kernels. Since $\spin$ is finite, \cvinf is compact, and a density
argument gives the conclusion. \qed
\end{pf}


\section{Applications\label{Applications}}
Using the increasing coupling, we develop a precise analysis
of the structure of the set of PCA dynamics' stationary measures. Moreover, the time-asymptotical behaviour
is investigated. \\ \indent
Let us first prove a property (Proposition~\ref{Sub_Super_Gibbs}) for
stationary measures associated to PCA restricted on $\cvinfLx{\L}$ ($\L\Subset \Zd$)
(see formula~(\ref{Def_PCA_fini})). %
The relation between these measures and the stationary measures for the PCA dynamics on \cvinf
is then established (Proposition~\ref{ResultatsDeMonotonie}).
In particular, we formulate an identity relating spatial limits and temporal limits
(see equations~(\ref{LimTempLimSpatPlus}) and~(\ref{LimTempLimSpatMinus})).
Proposition~\ref{Prop_VolFinVolInf} establishes a comparison between infinite volume PCA's ({\it e.g.} $\L=\mathbb Z^d$)
and PCA's in a large but finite volume.
 See also~\cite{PylExpErg} for applications. \\ \indent
In the following $P$ is an attractive PCA dynamics on \cvinf, where \spin is %
 a totally ordered space.

\subsection{Finite volume PCA dynamics\label{Sous_section_finite_volume}}
Let $\L \Subset \Zd$ be a finite subset of $\Zd$, called finite volume.
We call \impt{finite volume PCA dynamics with boundary
condition} $\tau$ ($\tau \in \cvinf$ or $\tau \in
\cvinfLx{\L^c}$),
the Markov Chain on
$\cvinfLx{\L}$ whose transition probability $P_\L^\t$
is defined by:
\begin{equation} \label{Def_PCA_fini}
P_\L^\t( d\s_\L \ | \ \e_\L\ )=\miniop{}{\otimes}{k \in \L} %
p_k(\ d\s_k \ | \ \e_\L\t_{\L^c} \ ) .
\end{equation}
It may be identified with the following infinite volume %
PCA dynamics on~$\cvinf$:
\be{DefinitionPCAVolFini}
 P_\L^\t( d\s \ | \ \e_\L\ )=\miniop{}{\otimes}{k \in \L} %
p_k(\ d\s_k \ | \ \e_\L\t_{\L^c} \ ) \otimes  \d_{\t_{\L^c }}
(d\s_{\L^c})
\ee
where the spins of $\L$ evolve according to $P_\L^\t$, and those of
$\L^c$ are almost surely `frozen' at the value~$\t$.
We assume that the finite volume PCA dynamics~$P_\L^\t$ are irreducible and aperiodic
Markov Chains. They then
admit one (and only one) stationary probability measure, called~$\nu_\L^\t$ %
(\ie $\nu_\L^\t P_\L^\t=\nu_\L^\t$); %
furthermore
$P_\L^\t$ is ergodic, which means $\lim_{n \to \infty} \rho_\L (P_\L^\t)^{(n)}= \nu_\L^\t$ %
in the weak sense, for any initial
condition~$\rho_\L$.\\ \indent
A sufficient condition for the irreducibility and aperiodicity of $P_\L^\t$ is for instance to assume that
the PCA dynamics under study are \impt{non degenerate}. This means:
$ \forall k \in \Zd,\ \forall \e \in \cvinf, \ \forall s \in \spin, \ %
p_k(\ s \ | \ \e \ ) >0 $.
The following Proposition states that %
the finite volume stationary measures associated with extremal boundary conditions
satisfy some sub/super-DLR relation (which means these measures are `sub/super-Gibbs measures').
In the very special case where $\spin=\{-1,+1\}$ and for $P$ reversible, this result was shown in~\cite{DPLR}.
\begin{prop} \label{Sub_Super_Gibbs}
Let $\nu_\L^\maximal$ (resp. $\nu_\L^\minimal$) be the unique stationary probability measure associated with
the finite volume PCA dynamics $P_\L^\maximal$ (resp. $P_\L^\minimal$) with %
$\maximalConf$ (resp. $\minimalConf$) extremal boundary condition.
Let \mbox{$\L \subset \L' \Subset \Zd$}. One has, for any configuration~$\s$, 
\begin{equation}
\nu_{\L'}^\minimal ( .| \s_{\L' \setminus \L} ) \succeq \nu_\L^\minimal (.)  \quad  \textrm{ and } \quad %
\nu_{\L'}^\maximal ( .| \s_{\L' \setminus \L} ) \preceq \nu_\L^\maximal (.) .
\end{equation}
\end{prop}
\begin{pf}
First, using~(\ref{Monotonie_locale}) shows that the pair of PCA's %
$(P_{\L'}^\maximal,P_\L^\maximal \otimes \delta_{\maximal_{\L' \setminus \L}})$
(resp. $(P_\L^\minimal \otimes \delta_{\minimal_{\L' \setminus \L}},P_{\L'}^\minimal)$ on $\cvinfLx{\L'}$)
is increasing.
Using the increasing coupling defined in Theorem~\ref{ExistenceCouplageCroissant}, we may then state: for any initial %
condition $\s$ and
for $n \geq 1$ that
$$ P_{\L'}^\maximal \opreceq \Big( P_\L^\maximal \otimes \delta_{\maximal_{\L' \setminus \L}} \Big) %
\big( f(\o^2(n)) - f(\o^1(n)) | \ %
(\o^1,\o^2)(0)=(\s,\s) \ \big) \geq 0,$$
where $f$ is any increasing function on \cvinf.
Thus
$$ P_{\L'}^\maximal  (f(\o(n)) \ | \ \o(0)=\s ) \leq %
   P_\L^\maximal \otimes \delta_{\maximal_{\L' \setminus \L}} (f(\o(n)) \ | \ \o(0)=\s ).$$
\noindent Letting $n \to \infty$ and
using finite volume ergodicity yields  %
$\nu_{\L'}^\maximal \preceq \nu_\L^\maximal \otimes \delta_{\maximal_{\L' \setminus \L}}$.
Similarly, $\nu_\L^\minimal \otimes \delta_{\minimal_{\L' \setminus \L}} \preceq \nu_{\L'}^\minimal $. \\ \indent
Let $\s_{\L' \setminus \L} \in \cvinfLx{\L' \setminus \L}$. Let $B$ be the event
$ B=\{ \o \in \cvinfLx{\L'} : \o_{\L' \setminus \L}= \s_{\L' \setminus \L} \}$.
Consider a sequence of independent, identically distributed random %
variables $(Z_n)_{n \geq 1}$, with distribution
$\nu_{\L'}^\maximal$. Let $Y$ be a random variable with distribution %
$\nu_\L^\maximal \otimes \delta_{\maximal_{\L' \setminus \L}}$.
Let $T$ be the stopping time $\inf \{ n \geq 1: Z_n \in B \}$. %
Checking that $ \forall  n \geq 1$,
$Z_n \preceq Y$ almost surely, one then has $Z_T \preceq Y$.
This in turn means that
$\nu_{\L'}^\maximal ( .| \s_{\L' \setminus \L} ) \preceq \nu_\L^\maximal (.)$ and
the other inequality is proved in the same way. \qed
\end{pf}
\begin{prop} \label{ResultatsDeMonotonie}
Let $\L \Subset \Zd$.
The measure
 $\nu_\L^{\maximalConf}$ (resp. $\nu_\L^{\minimalConf}$)
is  the maximal (resp. minimal) measure of the set
\mbox{$\{ \nu_\L^\t \tq \t \in \cvinfLx{\L^c} \}$}.
Let $\nusup$ and $\nuinf$ denote the maximal and the minimal elements %
of the
set~$\S$ of stationary measures on \cvinf associated to the PCA dynamics~$P$. \\ \indent
The following relations hold:
\be{LimTempLimSpatPlus}
\nusup= \lim_{L \to  \infty} \nu_\BL^{\maximalConf} \otimes %
\delta_{(\maximalConf)_{\BL^c}}=\lim_{n \to \infty} \delta_{\maximalConf} P^{(n)}
\ee %
\be{LimTempLimSpatMinus} %
\nuinf= \lim_{L \to  \infty} \nu_\BL^{\minimalConf} \otimes %
\delta_{(\minimalConf)_{\BL^c}}=\lim_{n \to \infty} \delta_{\minimalConf} P^{(n)},
\ee
where for $L$ integer, $\BL$ is the $l^1$-ball
$\{ k \in \Zd : \norme{k} = \sum_{j=1}^{d} | k_j | \leq 1 \}$.
In particular, $P$ admits a unique stationary measure~$\nu$ if and
only if %
\mbox{$ \nuinf=\nusup $}.
\end{prop}
\begin{pf}
Let us first prove that: \mbox{$ \t \preceq \t' \Rightarrow \nu_\L^\t \preceq \nu_\L^{\t'}$}.
Let $f$ be an increasing function on
$\cvinf$.
It is easy to check that $(P_\L^\t,P_\L^{\t'})$ is a increasing
pair, thus $P_\L^\t \opreceq P_\L^{\t'}$ preserves stochastic order.
Let $\s \in \cvinf$ be an initial condition. Since
\mbox{$ \s_\L\t_{\L^c} \preceq  \s_\L\t_{\L^c}'$},
such an inequality is  at time $n$ still valid. Using the monotonicity of~$f$, we have:
$$ P_\L^\t \opreceq P_\L^{\t'} %
\big( f(\o^2(n)) - f(\o^1(n)) | \ %
(\o^1,\o^2)(0)=(\s,\s) \ \big) \geq 0 \ . $$
Thus
$ P_\L^{\t} (f(\o(n)) \ | \ \o(0)=\s ) \leq %
   P_\L^{\t'} (f(\o(n)) \ | \ \o(0)=\s )  $.
The first result thus follows by letting $n \to \infty$ and
using  finite volume ergodicity; the extremality of $\nu_\L^\maximal$ and $\nu_\L^\minimal$ follows.\\ \indent
Then, note that
\mbox{$\lim_{L \to \infty} ( \nu_\BL^{\minimalConf} \otimes \delta_{(\minimalConf)_{\BL^c}})$} and %
\mbox{$\lim_{L \to \infty} (\nu_\BL^{\maximalConf} \otimes \delta_{(\maximalConf)_{\BL^c}})$} %
exist due to monotonicity of the following sequences:
\mbox{$ (\nu_\BL^- \otimes \delta_{(\minimalConf)_{\BL^c}})_{L}$} and %
 \mbox{$( \nu_\BL^+ \otimes \delta_{(\maximalConf)_{\BL^c}})_{L}$}.
This comes from the fact that
\mbox{$ \proj{\nu_{\L'}^{\maximalConf} }{\L} %
\preceq \nu_\L^{\maximalConf}  $} %
(where  $\L \Subset \L' \Subset \Zd$ and %
\mbox{$ \proj{}{\L}$} denotes  the projection
on $\L$) which is easily checked using the increasing coupling
\mbox{$(P_{\L'}^{\maximalConf},P_\L^{\maximalConf})$}.
Since $\nu_\BL^{\maximalConf}$ is $P_\L^{\maximalConf}$-stationary %
(resp. $\nu_\BL^{\minimalConf}$ is $P_\L^{\minimalConf}$-stationary),
the limits
\mbox{$\lim_{L \to \infty} ( \nu_\BL^{\minimalConf} \otimes \delta_{(\minimalConf)_{\BL^c}})$} and  %
\mbox{$\lim_{L \to \infty} (\nu_\BL^{\maximalConf} \otimes
  \delta_{(\maximalConf)_{\BL^c}})$} %
are $P$-stationary. \\ \indent
Let $\nu$ be a $P$-stationary measure, and $L$ any positive integer.
Since  the coupling \mbox{$P_\BL^- \opreceq P \opreceq P_\BL^+$}
preserves stochastic order, using finite volume ergodicity, one can state:\\
$\nu_\BL^{\minimalConf} \otimes \delta_{(\minimalConf)_{\BL^c}} \preceq \nu %
\preceq \nu_\BL^{\maximalConf} \otimes \delta_{(\maximalConf)_{\BL^c}}$.
We then have:
\be{IntermedEncadrementSpatial}
\lim_{L   \to \infty} \nu_\BL^{\minimalConf} %
\otimes \delta_{(\minimalConf)_{\BL^c}}  %
\preceq %
\nu \preceq \lim_{L   \to \infty} \nu_\BL^{\maximalConf} \otimes \delta_{(\maximalConf)_{\BL^c}}.
\ee
%
On the other hand, it is easy to check \mbox{$\delta_\maximalConf P \preceq
  \delta_\maximalConf$},
so that, using
$P$'s attractivity,
\mbox{$(\delta_\maximalConf P^{(n)})_{n \in \N}$} is
  decreasing. Analogously,
\mbox{$(\delta_\minimalConf P^{(n)})_{n \in \N}$} is increasing. Thus,
  the limits \mbox{$\lim_{n \to \infty} \delta_{\minimalConf} P^{(n)}$} and  %
\mbox{$\lim_{n \to \infty} \delta_{\maximalConf} P^{(n)}$} exist and
are obviously $P$-stationary measures. \\ \indent
%
Let $\nu$ be a $P$-stationary measure. Since $P$ is
  attractive and \mbox{$\delta_{\minimalConf} %
 \preceq \nu \preceq \delta_{\maximalConf} $}, we have:
 \be{EncadrementTemporel}
\lim_{n \to \infty} \delta_{\minimalConf} P^{(n)} %
\preceq \nu \preceq \lim_{n \to \infty} \delta_{\maximalConf} P^{(n)}.
\ee
Using the fact that the measures
 \mbox{$\lim_{L \to \infty} ( \nu_\BL^{\minimalConf} \otimes \delta_{(\minimalConf)_{\BL^c}})$}, %
\mbox{$\lim_{L \to \infty} (\nu_\BL^{\maximalConf} \otimes
  \delta_{(\maximalConf)_{\L^c}})$},  
\mbox{$\lim_{n \to \infty} \delta_{\minimalConf} P^{(n)}$} and
\mbox{$\lim_{n \to \infty} \delta_{\maximalConf} P^{(n)}$} are
$P$-stationary,
we apply inequalities~(\ref{IntermedEncadrementSpatial})
and~(\ref{EncadrementTemporel}), and the
conclusion follows. \qed
\end{pf}
%
%
%
\subsection{Comparison of finite \& infinite volume PCA\label{Sous_section_fini_infini}}
Thanks to the above constructed coupling, we investigate the time-asymptotical behaviour.
The PCA dynamics $P$ on the infinite volume space $\cvinf$ considered in this subsection
is assumed to be translation invariant (or \impt{space homogeneous}):
\mbox{$\forall k \in \Zd$}, \mbox{$\forall s \in \spin$}, \mbox{$\forall \e \in \cvinf$}, %
$p_k(\ s \ | \ \e \ )= %
p_0(\ s \ | \ \theta_{-k}\e \ )$,
where \mbox{$\theta_{k_0}(\s) =$} \mbox{$ (\s_{k-k_0})_{k \in \Zd}$}.
Remark that %
if
the PCA dynamics $P^i$ are translation invariant, so is the coupled dynamics $P^1 \opreceq \hdots \opreceq P^N$. \\
\indent We will use
the notation $\coup{}$
to denote the coupling $P \opreceq P \opreceq \hdots \opreceq P$
of $N$~times the same attractive PCA dynamics~$P$.
Using the compatibility property of the Lemma~\ref{Propriete_Compatibilite},
the marginal
of $P^{\opreceq N'}$
on $N$ components chosen in \mbox{$\{1,..,N'\}$} is the same as the coupling
$P^{\opreceq N}$.
Here, it is enough to choose $N=4$. \\
\indent As Proposition~\ref{ResultatsDeMonotonie} shows, in order to study
the behaviour of a PCA dynamics $P$ on $\cvinf$, one may turn its attention %
to the finite volume associated dynamics $P_\L^\t$ on $\cvinfLx{\L}$, where $\L \Subset \Zd$.
Note that their time asymptotics are known.
Using the increasing coupling, Proposition~\ref{Prop_VolFinVolInf} below
shows how the time-asymptotical behaviour of our PCA
is controlled by the sequence $(\rho(n))_{n \geq 1}$,
where
\be{DefRho}
\r(n)=\coupling{\o_0^1(n) \neq \o_0^2(n)}{(\o^1,\o^2)(0)%
=(\minimalConf,\maximalConf)},
\ee
where $\coup{}$ is the coupling introduced in Corollary~\ref{Corollary_main}. In
the paper~\cite{PylExpErg} we gave conditions to ensure the convergence of
$(\rho(n))_{n \geq 1}$ and stated conditions for the ergodicity with exponential speed of the dynamics~$P$. \\
\indent Let $\L \Subset \Zd$.
Let $P_\L^\maximalConf$ (resp. $P_\L^\minimalConf$)
be
the dynamics on $\cvinfLx{\L}$ defined in~(\ref{DefinitionPCAVolFini})
with the maximal (resp. minimal) boundary condition $\maximalConf$
(resp.~$\minimalConf$). First note the easily checked fact:
\begin{lem}
If the PCA dynamics $P$ is attractive then
$(P_\L^\minimalConf,P,\hdots,P,P_\L^\maximalConf)$
is increasing, and thus the increasing coupling
$ P_\L^\minimalConf \opreceq P \opreceq \hdots \opreceq P \opreceq P_\L^\maximalConf$
can be defined.
\end{lem}
\begin{prop} \label{Prop_VolFinVolInf}
Let $\s \preceq \e \in \cvinf$ and $P$ be an attractive PCA dynamics.
The following
inequality holds:
\begin{equation} \label{encadrement}
\coupling{\o_0^1(n) \neq \o_0^2(n)}{(\o^1,\o^2)(0)=(\s,\e)} %
 \leq \r(n) \leq P_\L^\minimalConf \opreceq P_\L^\maximalConf %
 ( \o^1_0(n)  \neq \o^2_0(n) \ | %
(\o^1,\o^2)(0)=(\minimalConf,\maximalConf) )
\end{equation}
where
$(\r(n))_{n \in \mathbb N^*}$ is defined by~(\ref{DefRho}).
For each initial condition $\xi$ on $\cvinf$ and for any time
 $n$, it holds:
\begin{equation}  \label{PCAloisOrdonnees}
P_\L^{\minimalConf} \big( \o (n) \in .  \big| %
\o (0)=\xi_\L ({\minimalConf})_{\L^c} \big) %
 \preceq %
P \big( \o(n)\in .  \big  |   \o(0)=\xi \big) \preceq  %
 P_\L^{\maximalConf} \big( \o(n) \in .  \big|  \o(0)=\xi_\L %
({\maximalConf})_{\L^c} \big) .%
\end{equation}
The sequence
$(\r(n))_{n \in \mathbb N^*}$ is decreasing, and
$P$ is ergodic if and only if \mbox{$\lim_{n \to \infty} \r(n)=0$}.
Moreover, in this case,
\begin{equation} \label{Ineq_rho_ergodicite}
 \sup_{\s} \abs{ %
P\Big( f(\o ( n )) | \o(0)=\s \Big) - \nu(f) }%
  \leq %
 2 \ \ntrois{f} \ \r(n)
\end{equation}
where $\nu$ denotes the unique stationary measure and
where, for each $f$ continuous function on the compact \cvinf and
for all $k$ in $\Zd$,  %
$\Var{k}{f}=\sup \Big\{ \abs{f(\s)-f(\e)} \ : \ (\s,\e)\in (\cvinf)^2, %
\s_{\{ k  \}^c} \equiv \e_{\{ k  \}^c} \Big\} $,
whereas  %
\mbox{$ \ntrois{f} = \sum_{k \in \Zd} {\Var{k}{f}}$}.
\end{prop}
\begin{pf}
The proof of the left inequality in~(\ref{encadrement}) %
 is straightforward using  the compatibility property from Lemma~\ref{Propriete_Compatibilite}.
The right inequality comes from the confunction of two properties: preservation of the stochastic order as well as
the %
compatibility property of the coupling %
$ P_\L^{\minimalConf} \opreceq P \opreceq P \opreceq P_\L^{\maximalConf}$.\\ \indent
Since  the coupling
$P_\L^- \opreceq P \opreceq P_\L^+$ is increasing,
(\ref{PCAloisOrdonnees}) is a consequence of the fact that
any initial condition  $\xi$ in $\cvinf$ satisfies
\mbox{$ \xi_\L(\minimalConf)_{\L^c}  \preceq \xi %
\preceq  \xi_\L (\maximalConf)_{\L^c} $}.\\ \indent
The monotonicity of the sequence $(\r(n))_{n \in \mathbb N^*}$
comes from the coalescence property of the increasing coupling
\mbox{$\coup{}$}. \\ \indent
If $P$ is ergodic, there can be only one stationary measure on \cvinf and so $\lim_{n \to \infty} \r(n)=0$. \\ \indent
Conversely, let $f$ be a local function.
For any $\s,\e$ configurations in \cvinf, let us write:
\begin{eqnarray}
\lefteqn{   \abs{ %
P(f(\o(n)) | \o(0)=\s ) - P(f(\o(n)) | \o(0)=\e ) } } \nonumber \\
& \leq & %
\abs{ \coup{} \Big( f(\o^1(n)) -f(\o^2(n)) \ \Big | %
(\o^1,\o^2)(0)=(\minimalConf,\s) \Big) }  \label{IntermedDeux} %
 + \abs{ \coup{} \Big( f(\o^1(n)) -f(\o^2(n)) \ \Big | %
(\o^1,\o^2)(0)=(\minimalConf,\e) \Big) } \nonumber .
\end{eqnarray}
Since $f$ is local, for all $\xi^1,\xi^2$,
\mbox{$\abs{f(\xi^1)-f(\xi^2)}$}
depends only on $\xi^1_{\L_f}$ and $\xi^2_{\L_f}$,
which differ only in a finite
number of sites.
Using interpolating configurations between
 $\xi^1_{\L_f}$ and $\xi^2_{\L_f}$, we write: \\
\mbox{$|{f(\xi^1)-f(\xi^2)}| \leq \sum_{k \in \L_f} \Var{k}{f} %
\ind{\{\s_k\neq  \e_k\}}{} $}, so that the translation invariance assumption and
the left part of~(\ref{encadrement}) then yield:
$\abs{ %
P(f(\o(n)) | \o(0)=\s ) - P(f(\o(n)) | \o(0)=\e ) } %
 \leq \  2 \ \ntrois{f} \r(n)$,
which is enough to conclude and state~(\ref{Ineq_rho_ergodicite}). \qed
\end{pf}
%
%

\section{Partially ordered spin space case\label{section_ordre_partiel}}
In all this section, \spin is a partially ordered space.
When $\spin$ is totally ordered, a necessary and sufficient condition for the existence of an increasing
coupling of PCA dynamics is given in section~\ref{proofs_main} by the inequality~(\ref{ordre_fonction_repart}).
It is done in term of
the distribution function $F_k(.,\s)$ ($\s$ given) of the probability $p_k(\ . \ | \s)$. %
We recall previous results,
 who gave a
necessary and sufficient condition
for the existence of an increasing
coupling of two PCA's. In particular, $P$ is attractive if and only if
\begin{equation} \label{Condition_simple_Lopez_Sanz}
\forall k \in \Zd, \forall \s \preceq \e, \forall \ \Gamma \textrm{ up-set in } \spin \Rightarrow %
\miniop{}{\sum}{s \in \Gamma} p_k(s|\s) \leq  \miniop{}{\sum}{s \in \Gamma} p_k(s|\e) .
\end{equation}
The quantity which now makes sense is the generalised
function $F_k(\Gamma,\s):= \sum_{s' \in \Gamma} p_k(s'|\s)$, where $\Gamma$ is an up-set of $\spin$
(see Definition~\ref{def_upset}). \\
\indent Nevertheless, there is a gap between coupling two PCA's or more than three PCA's.
The counter-examples A and B presented here show that %
a satisfactory coupling of three PCA's
may not exist and
condition~(\ref{Condition_simple_Lopez_Sanz}) of \cite{LopezSanz2000}
is not
sufficient
for the existence of an increasing $3$-coupling when $\spin$ is any partially ordered space. \\ \indent
These counter-examples rely on examples 1.1 and 5.7 in \cite{FillMachida}
of stochastically monotone families of distributions, indexed by a partially ordered set, which are not
realisable monotone in the following sense.
Let $(Q_\alpha)_{\alpha \in \Aposet}$ be a family of probability
distributions on a finite set $\Sposet$
indexed by a partially ordered set $\Aposet$.
\cite{FillMachida} define the system $(Q_\alpha)_{\alpha \in \Aposet}$ as
stochastically monotone if $\alpha_1 \Apreceq \alpha_2$ implies $Q_{\alpha_1} \Spreceq Q_{\alpha_2}$.
It is said to be realisable monotone if there exists a system of $\spin$-valued random variables %
$(X_\alpha)_{\alpha \in \Aposet}$ , defined
on the same probability space, such that the distribution of $X_\alpha$ is $Q_\alpha$ and
$\alpha_1 \Apreceq \alpha_2$ implies $X_{\alpha_1} \Spreceq X_{\alpha_2}$ a.s. \\ \indent
In our case, the existence of a coupling of the $N$ PCA dynamics $(P^1,\hdots,P^N)$
implies that, for any $k \in \Zd$ fixed, the system $\{ p_k(.|\s_{V_k}) : {\s_{V_k} \in \cvinfLx{V_k}} \}$ of probability distributions on $\spin$,
which is indexed by the partially ordered set
$\cvinfLx{V_k}$, is realisable monotone.
In the counter-examples presented here, the distributions %
are stochastically monotone but not realisable monotone. \\
\begin{defn}\label{def_upset}
A subset $\Gamma$ of \spin is said to be an \impt{up-set} (or increasing set) (resp. down-set or decreasing set) %
 if: $x \in \Gamma, y \in \spin, x \preceq y \Rightarrow y \in \Gamma$ (resp. %
 $x \in \Gamma, y \in \spin, x \succeq y \Rightarrow y \in \Gamma$).
\end{defn}
\indent Note that the indicator function of an up-set (resp. down-set) is an increasing (resp. decreasing)
function. Moreover, Theorem~1 in \cite{KKOB} states that for two measures $\mu_1,\mu_2$ on \spin, %
 $\mu_1 \preceq \mu_2$ if and only if
$\mu_1(\Gamma) \leq \mu_2(\Gamma)$ for all up-sets $\Gamma$ of $\spin$, which is equivalent
to $\mu_1(\Gamma) \geq \mu_2(\Gamma)$ for all down-sets $\Gamma$ of $\spin$.

\begin{exempleNUM} \label{diamant}
\indent Let $\spin=\spin_A=\{0,1 \}^2$ be equipped with the natural partial order
$(0,0) \preceq (0,1)$, $(0,0) \preceq (1,0)$,
$(0,1) \preceq (1,1)$, $(1,0) \preceq (1,1)$
(where $(0,1)$ and $(1,0)$ are not comparable).
Let $P=\miniop{}{\otimes}{k \in \Zd} p_k$, %
where $p_k(\ . \ | \s) =p_k(\ . \ | \s_k)$ is defined as follows:
\begin{equation} \label{contre_exemple_1}
\left\{ \begin{array}{l}
p_k(\ . \ | (0,0))= \frac{1}{2} (\delta_{(0,0)}+ \delta_{(1,0)}) \\
p_k(\ . \ | (1,0))= \frac{1}{2} (\delta_{(0,0)}+ \delta_{(1,1)}) \\
p_k(\ . \ | (0,1))= \frac{1}{2} (\delta_{(0,1)}+ \delta_{(1,0)}) \\
p_k(\ . \ | (1,1))= \frac{1}{2} (\delta_{(1,0)}+ \delta_{(1,1)}) . \\
\end{array} \right.
\end{equation}
It is simple to check that this PCA dynamics is attractive. Nevertheless, $P^{\opreceq 4}$ can not exist since
$(p_k(\ . \ | (0,0)), p_k(\ . \ | (1,0)), %
p_k(\ . \ |(0,1)), p_k(\ . \ | (1,1)))$ is a stochastically monotone family which is not
realisable monotone (see example~1.1 in \cite{FillMachida}).
This PCA is in fact a collection of independent, $\spin$-valued Markov Chains,
whose transition probability is
$p_0(.|.)$.
This example states the non-existence of a coupling of four particular Markov Chains.
\end{exempleNUM}
\begin{exempleNUM} \label{ordre_Y}
\indent Let $\spin=\spin_B=\{x,y,z,w \}$, considered with the following partial order
$x \preceq z, y \preceq z, z \preceq w$ and $x$ and $y$ are not comparable.
Consider the dimension $d=1$, and the PCA $P=\miniop{}{\otimes}{k \in \mathbb Z} p_k$ %
where $p_k(\ . \ | \s) =p_k(\ . \ | \s_{\{k,k+1 \}})$ is defined as follows:
\begin{equation}
\left\{ \begin{array}{lcl} \label{contre_exemple_2}
p_k(\ . \ | (x,y))= \frac{1}{2} (\delta_{x}+ \delta_{y}) %
& \qquad & p_k(\ . \ | (y,z))= \delta_{z} \\
p_k(\ . \ | (x,z))= \frac{1}{2} (\delta_{x}+ \delta_{w}) %
 & \qquad & p_k(\ . \ | (z,x))= \delta_{z} \\
p_k(\ . \ | (z,y))= \frac{1}{2} (\delta_{y}+ \delta_{w}) %
 & \qquad & p_k(\ . \ | (x,x))= \delta_{x} \\
p_k(\ . \ | (z,z))= \frac{1}{2} (\delta_{z}+ \delta_{w}) %
& \qquad & p_k(\ . \ | (y,x))= \delta_{z} \\
p_k(\ . \ | (y,y))= \delta_{y} %
& \qquad & p_k(\ . \ | \textrm{otherwise})= \delta_{w} \\
\end{array} \right.
\end{equation}
This is an attractive PCA, nevertheless a synchronous coupling $P^{\opreceq 4}$ can not exist since \\
$(p_k(\ . \ | (x,y)), p_k(\ . \ | (x,z)), %
p_k(\ . \ | (z,y)), p_k(\ . \ | (z,z)))$
is a stochastically monotone family which is also non realisable monotone %
(see example~5.7 in \cite{FillMachida}).
\end{exempleNUM}

Let us now present some generalisation of our main results, %
Theorem~\ref{ExistenceCouplageCroissant} and
Corollary~\ref{Corollary_main}, when the spin space %
\spin\ belongs to a special class~$\mathcal Z$ of partially ordered sets introduced %
in~\cite{FillMachida} and called \impt{linearly ordered spaces}. \\
\indent We call \impt{predecessor} (resp. \impt{successor}) of $s$ ($s \in \spin$) %
any element $s'$ such that $s \preceq s'$ %
(resp. $s \succeq s'$) and
$s \preceq s'' \preceq s' \Rightarrow$ $s'' \in \{s,s'\}$ (resp. %
$s \succeq s'' \succeq s' \Rightarrow$ $s'' \in \{s,s'\}$).
$\spin$ belongs to the class~$\mathcal Z$ if, for any $s \in \spin$, only one of the following situations occurs:
$s$~admits exactly one successor and one predecessor;
 $s$~admits no predecessor and at most two successors;
$s$~admits no sucessor and at most two predecessors.
One can then define on $\spin$ a linear order~$\preceqlino$ by %
numbering the elements of $\spin$:
$\{s_1,\hdots,s_n\}$ (where $n=\card{\spin}$) in such a way that
$s_{i+1}$ be a sucessor or a predecessor of $S_i$ (for $i=1,\hdots,n$) and
declaring that $s_i \preceqlino s_j$ if $i \leq j$.
Of course such linear order might be incompatible with the partial order $\preceq$ originally %
defined on $\spin$. \\
The set $\spin_C=\{s_i, 1 \leq i\leq 10 \}$ with order relations:
$s_1 \preceq s_2 \preceq s_3 \preceq s_4$, $s_6 \preceq s_5 \preceq s_4$,
$s_6 \preceq s_7 \preceq s_8$, and $s_{10} \preceq s_9 \preceq s_8$ is an example of such a space.
On the other hand, the spin
space $\spin_D=\{x,y,z,u,v,w\}$ with order relations
$ y\preceq z$, $x\preceq z$, $w \preceq z$, $w \preceq u$, $w \preceq v$
does not belong to this class. \\
\indent
Define, for $s_i \in \spin$ ($1 \leq i \leq n$), the subset \semi{s_i} of \spin\ with
$\semi{s_i}=\{ s_j \in \spin \ : s_j \preceqlino s_i \}.$
The sets $\semi{s}$ (with $s \in \spin$) are either up-sets or
down-sets of $\spin$ (respective to the original order $\preceq$ on \spin). For instance, when $\spin=\spin_C$, $\semi{s_5}$ is an up-set and
$\semi{s_6}$ is a down-set.
It suffices to consider the generalised function $F_k(\Gamma,\s)$ for sets
of the form $\Gamma=\semi{s}$ for which we have:
\begin{equation}
F_k(s,\s )= p_k(\ \semi{s} \ | \ \s )= %
\miniop{}{\sum}{s' \in \semi{s}} p_k(s' | \s) \qquad (s \in \spin, \ \s \in \cvinf).
\end{equation}
\indent When $\spin$ is a linearly ordered set, the %
monotonicity condition~(\ref{Monotonie_locale}) %
is equivalent to the following conditions for the generalised associated
distribution functions (Lemma~5.5 in~\cite{FillMachida}): \\
$  \forall k \in \Zd, \forall (\z^1,%
\z^2,\hdots,\z^N) \in (\cvinf)^N \textrm{ such that } %
\ \z^1 \preceq \z^2 \preceq \hdots \preceq \z^N$,
\begin{gather}
\label{cond_monoton_class_Z_1}
\forall s \in \spin \textrm{ with } \semi{s} \textrm{ down-set, }
F_k^1(s \ | \ \z^1 ) \geq F_k^2(s \ | \ \z^2 ) \geq \hdots %
\geq F_k^N(s \ | \ \z^N ). \\
%
\label{cond_monoton_class_Z_2}
\forall s \in \spin \textrm{ with } \semi{s} \textrm{ up-set, }
F_k^1(s \ | \ \z^1 ) \leq F_k^2(s \ | \ \z^2 ) \leq \hdots %
\leq F_k^N(s \ | \ \z^N )
\end{gather}
\begin{prop} %
When \spin is a linearly ordered spin space, Theorem~\ref{ExistenceCouplageCroissant} and
Corollary~\ref{Corollary_main} still hold. \end{prop}
\begin{pf} The proof of such results relies on the following construction.
Let us define the generalised probability transform, for $\s \in \cvinf$ and $k\in \Zd$ fixed:
$ (F_k(\ . \ ,\s))^{-1} (t) = \miniop{}{\inf}{\preceqlino} \{ s_k : t < F_k(s_k,\s) \}$, ($t \in ]0,1[$),
where the infimum is given in term of the linear order~$\preceqlino$.
The construction of the increasing coupling holds as before thanks to the following %
evolution rule~(\ref{regle_locale_algo}) between times~$n$ and~$n+1$, where $(F_k^i)^{-1}$ denotes the generalised
distribution function, as introduced before.
The coherence of this coupling %
with the partial order $\preceq$ is proved in Lemma~6.2 in~\cite{FillMachida}. \qed
\end{pf}

\begin{ack}
The author thanks P.~Dai~Pra and
G.~Posta for helpful comments, O.~Häggström for mentioning the reference~\cite{FillMachida}, and
S.~R\oe lly for a careful reading and comments to the preliminary version of this paper.
He thanks as well M.~Sortais for helping to improve the english exposition.
\end{ack}
\newpage
\bibliographystyle{elsart-harv}

\end{document}